\documentclass{birkjour}

\usepackage{graphicx}
\usepackage{amssymb}
\usepackage{amsthm}

\theoremstyle{remark}

\newcommand{\perim}{{\rm perim}}

\begin{document}

\baselineskip= 19.5pt

\title[Three characterizations of reducedness of spherical convex polygons]{Three characterizations of reducedness of spherical convex polygons}

\author[M. Lassak]{Marek Lassak}

\address
{University of Science and Technology\\
85-789 Bydgoszcz, Poland}

\email{lassak@pbs.pl}

\date{}

\begin{abstract}
\noindent 
Denote by $S^2$ the two-dimensional sphere.
A spherical convex body on $S^2$ which does not properly contain a spherical convex body of the same spherical thickness is called a reduced body.
We give three characterizations of reducedness of spherical convex odd-gons on $S^2$.
Analogous characterizations hold true also in the Euclidean and hyperbolic planes.
\end{abstract}

\keywords{Sphere, reduced body, spherically convex body, spherically convex polygon, width, thickness}

\maketitle

{ }
\vskip-0.8cm

\textbf{MSC:} 52A55.

\section{Introduction}
\noindent
Denote by $S^2$ be the unit sphere of the three-dimensional Euclidean space $E^3$. 
By a {\it great circle} we mean the intersection of $S^2$ with any two-dimensional subspace of $E^{3}$. 
By a pair of {\it antipodes} we understand any pair of points obtained as the intersection of the sphere $S^2$ with a one-dimensional subspace of $E^3$.

If two different points $a, b$ are not antipodes, there is exactly one great circle containing them.
By the {\it spherical arc} $ab$, or shortly {\it arc} $ab$ we understand the shorter part of the great circle connecting the points $a$ and $b$. 
By the {\it distance} $|ab|$ of these points we mean the length of the arc connecting them. 

A subset of $S^2$ is called {\it convex} if it does not contain any pair of antipodes and if together with every two points it contains the arc connecting these points.
By {\it a convex body} on $S^2$ we understand a closed convex set whose interior is  non-empty. 

Any of the two subset of $S^2$ bounded by a great circle is called a {\it hemisphere}.
By its {\it center} we mean its point which is equidistant from all points of thie great circle.

If the centers of different hemispheres $G$ and $H$ are not antipodes, then $L = G \cap H$ is called a {\it lune}. 
The semicircles bounding $L$ and contained in the boundaries of $G$ and $H$ are denoted by $G/H$ and $H/G$, respectively. 
The {\it thickness} $\Delta (L)$ of $L$ is defined as the distance of the centers of $G/H$ and $H/G$.  

Recall a few notions and facts presented in \cite{L2} and repeated in the survey article \cite{L4}.
If a great circle $Q$ has a common point with a convex body $C \subset S^2$ but not with its interior, it is said that $Q$ {\it supports} $C$.
If $H$ is a hemisphere bounded by $Q$ and containing $C$, we say that it {\it supports} $C$.
By {\it the width ${\rm width}_K (C)$ of $C$ determined by} $K$ we mean the thickness of the lune of the form $K \cap K^*$ of the minimum thickness, where $K^*$ is a hemisphere supporting $C$. 
By compacness arguments we gather rhat such a $K^*$ exists.
What is more, from Part I of Theorem 1 from \cite{L2} we conclude that such $K^*$ is unique if the minimum is below $\frac{\pi}{2}$).
By the {\it thickness} $\Delta (C)$ of $C$ we understand the minimum width of $C$ determined by $K$ over all supporting hemispheres $K$ of $C$. 
The thickness of $C$ is nothing else but the minimum thickness of a lune containing $C$.

The convex hull $V$ of $k \geq 3$ points on $S^2$ such that each of them does not belong to the convex hull of the remaining points is called a {\it spherically convex $k$-gon}. 
These $k$ points are called {\it vertices} of $V$. 
We write $V= v_1v_2\dots v_k$ provided $v_1, v_2, \dots , v_k$ are successive vertices of $V$ when we are moving around $V$ on the boundary of $V$ according to the positive orientation. 
When the distances of every two successive vertices are equal, we obtain a {\it regular spherical $k$-gon}. 

After \cite{L2} we say that a spherical convex body $R \subset S^2$ is {\it reduced} provided $\Delta (Z) < \Delta (R)$ for every convex body $Z \subset R$ different from $R$ (this is an analog of the notion of reduced body in $E^2$ given by Heil in \cite{He}). 
The paper \cite{LM} is devoted to reduced bodies in $S^2$. 
Simple examples of reduced spherical convex bodies are spherical bodies of constant width and, in particular, the disks. 
There exists a large class of reduced odd-gons. 
From \cite{L3} we know that every regular spherical odd-gon of thickness at most $\pi \over 2$ is a reduced spherical body. 
  
Let $a$ be a point in a hemisphere different from its center and let $F$ be the great circle bounding the hemisphere. 
By {\it the projection of $a$ on} $F$ we mean this point $p \in F$ such that the distance $|ap|$ is the smallest from amongst all $|ac|$, where $c \in F$.

Two sets on $S^2$ are called to be {\it symmetric with respect to a great circle} provided in $E^3$ they are symmetric with respect to the plane of $E^3$ containing this circle.

The aim of this note is to give some characterizations of reducedness of convex polygons on $S^2$. 
Analogous characterizations are true also in the Euclidean and hyperbolic planes.

\section{Reduced spherical polygons}

After the paper \cite{L3} recall that every reduced spherical polygon is of thickness at most~$\pi \over 2$.

For a convex odd-gon $V = v_1v_2\dots v_n$, by the {\it opposite side} to the vertex $v_i$ we mean the side $v_{i + (n-1)/2}v_{i + (n+1)/2}$.
Here and later, always the indices are taken modulo~$n$.

In this note we deal with reduced spherical polygons of thickness below $\pi \over 2$.
By the way, spherical polygons of thickness $\pi \over 2$ are considered by Chang, Liu and Su \cite{CLS}, and also by Hue \cite{Ha}.

Recall Theorem 3.2 from \cite{L3}. 
{\it Every reduced spherical polygon is an odd-gon of thickness at most $\pi \over 2$.
A spherically convex odd-gon $V$ with $\Delta (V) < {\pi \over 2}$ is reduced if and only if the projection of every its vertices on the great circle containing the opposite side belongs to the relative interior of this side and the distance of this vertex from this side is $\Delta (V)$.}

The mentioned projection of $v_i$, where $i \in \{1, \dots, n\}$, onto the opposite side $v_{i + (n-1)/2}v_{i + (n+1)/2}$ is denoted by $t_i$.

The following figure from \cite{L2} which helps us to visualize the further explanations of this note.

\vskip0.4cm
\begin{center}
\includegraphics[width=7.45cm,height=7.45cm]{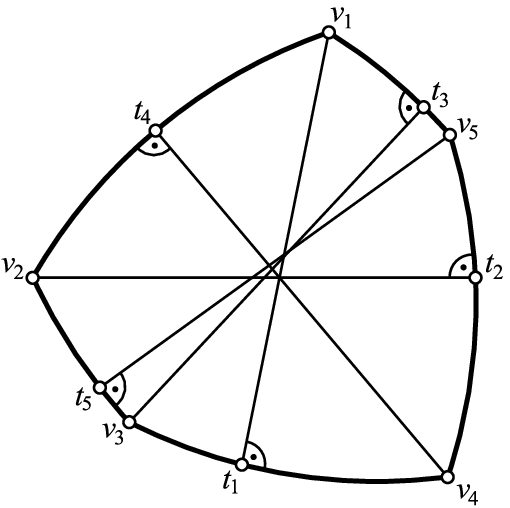}\\
\vskip0.3cm
\centerline
{Fig. A reduced spherical pentagon}
\end{center}
\vskip0.3cm

Inspired by the paper \cite{CLS} of Chang Liu and Su, the author put the following Problem 2.4 in \cite{L4}.
{\it Consider an arbitrary convex odd-gon $V= v_1v_2 \dots v_n$ on $S^2$ of thickness below $\frac{\pi}{2}$ which satisfies $|v_it_{i + (n+1)/2}| = |t_iv_{i + (n+1)/2}|$ for every $i \in \{1, \dots , n\}$.  
Is $V$ a reduced polygon?}

This problem rendered us to prove the following equivalence of (a) that $V$ is reduced to three other conditions.
In particular, the below implication (b) $\Rightarrow$ (a) gives a solution of this problem.

\vskip0.3cm
\noindent
{\bf Theorem.}
{\it On the sphere $S^2$ let us consider an arbitrary convex odd-gon $V= v_1v_2 \dots v_n$ of thickness below $\frac{\pi}{2}$ such that for every $i \in \{1, \dots, n\}$ the projection of $v_i$ on the straight line containing the side $v_{i + (n-1)/2}v_{i + (n+1)/2}$ belongs to the relative interior of this side.
The following conditions are equivalent:

{\rm (a)}\; $V$ is reduced,

{\rm (b)}\; $|v_it_{i + (n+1)/2}| = |t_iv_{i + (n+1)/2}|$ for $i = 1, \dots, n$,

{\rm (c)}\; $\angle t_{i + (n+1)/2}v_it_i = \angle t_iv_{i + (n+1)/2}t_{i + (n+1)/2}$ for $i = 1, \dots, n$,

{\rm (d)}\; the arc $v_it_i$ halves the perimeter of $V$ for $i = 1, \dots, n$.}

\begin{proof}
(a) $\Rightarrow$ (b), (c) and (d). 
This follows by Corollaries 3.6, 3.7 and 3.8 of~\cite{L3}.

\vskip0.15cm
(b) $\Rightarrow$ (a). 
Clearly, the arcs $v_it_i$ and $v_{i + (n+1)/2}t_{i + (n+1)/2}$ intersect for every $i \in \{1, \dots n\}$. 
By $o_i$ let us denote the point of intersection.
Since these arcs intersect, we have $\angle v_io_it_{i + (n+1)/2} = \angle t_io_iv_{i + (n+1)/2}$.
Of course, $\angle v_it_io_i = 90^\circ = \angle v_{i + (n+1)/2}t_io_i$.
Moreover, by our assumption we have $|v_it_{i + (n+1)/2}| = |t_iv_{i + (n+1)/2}|$.
Hence the triangles $v_it_{i + (n+1)/2}o_i$ and $t_iv_{i + (n+1)/2}o_i$ are congruent.
So $|v_io_i| = |v_{i + (n+1)/2}o_i|$ and $|t_{i + (n+1)/2}o_i| = |t_io_i|$.
As a consequence we obtain
$|v_io_i| + |o_it_i| = |v_{i + (n+1)/2}o_i| + |t_{i + (n+1)/2}o_i|$. 
Thus we get $\Delta_i = \Delta_{i + (n+1)/2}$, where $\Delta_i$ denotes the width of $V$ determined by this hemisphere supporting $V$ whose bounding semicircle contains the side $v_iv_{i+1}$.
Therefore also $\Delta_{i + (n+1)/2} = \Delta_{i + (n+1)/2 + (n+1)/2}$.
The last is nothing else but $\Delta_{i+1}$, which implies $\Delta_i = \Delta_{i+1}$ for every $i \in \{1, \dots n\}$. 
We conclude that $\Delta_1 = \dots = \Delta_n$. 
So $V$ is reduced. 

\vskip0.15cm
(c) $\Rightarrow$ (a).
Look to the proof of the implication (b) $\Rightarrow$ (a). 
Observe that the assumption $\angle t_{i + (n+1)/2}v_it_i = \angle t_iv_{i + (n+1)/2}t_{i + (n+1)/2}$ in place of the assumption $|v_it_{i + (n+1)/2}| = |t_iv_{i + (n+1)/2}|$ there is sufficient to conclude that the triangles $v_it_{i + (n+1)/2}o_i$ and $t_iv_{i + (n+1)/2}o_i$ are congruent.
The rest of the proof is as in the proof of the implication (b) $\Rightarrow$ (a). 

\vskip0.15cm
(d) $\Rightarrow$ (a). 
Consider the arcs $v_it_i$ and $v_{i+ (n+1)/2}t_{i+(n+1)/2}$.
Each of them halves the perimeter of $V$.
When going around the boundary of $V$ counterclockwise from  $v_i$ to $t_i$, the length is a half of perimeter.
The same when we go from $v_{i + (n+1)/2}$ up to $t_{i + (n-1)/2}$.  
These two pieces of the boundary differ at two arcs.
Namely, $v_it_{i + (n+1)/2}$ is covered two times, while the relative interior of $t_iv_{i + (n+)/2}$ is not covered at all.
Consequently, $\frac{1}{2}\perim (V) - |v_it_{i + (n+1)/2}| + \frac{1}{2}\perim (V) + |t_iv_{i + (n+1)/2}| = \perim (V)$.
Thus $|v_it_{i + (n+1)/2}| = |t_iv_{i + (n+1)/2}|$ for every $i \in \{1, \dots , n\}$.  
Hence (b) holds true.
By the implication (b) $\Rightarrow$ (a) the polygon $V$ is reduced. 
\end{proof}

Observe that Theorem remains true also for reduced polygons in the Euclidean plane $E^2$ and for reduced ordinary polygons in the hyperbolic plane $H^2$.
The proofs are analogous; we are using here strips instead of lunes.
For reduced polygons in $E^2$ see \cite{L1}, and for ordinary reduced polygons in $H^2$ see \cite{L5} and the paper \cite{S} by \'A. S\'agmeister.

\end{document}